\documentclass[11pt]{article}
\usepackage[T1]{fontenc}
\usepackage{amsmath}
\usepackage{amssymb}
\usepackage{graphicx}
\usepackage{subfig}
\usepackage{booktabs}
\usepackage{hyperref}
\usepackage{geometry}
\usepackage{booktabs}
\usepackage{tabularx}

\newcommand{\keywords}[1]{\textbf{Keywords:} #1}

\title{Sustainable and Optimal Harvesting in a Seasonally Harvested Fishery with a Marine Protected Area: A Two-Patch Model with Bang-Bang and Singular Control}

\author{Dinesh Kumar} % Replace with your name or details for arXiv submission

\date{}

\begin{document}

\maketitle

\begin{abstract}
We analyze a bioeconomic model for optimal fishery harvesting in a spatially heterogeneous habitat comprising both harvestable and preservation (reserve) zones. The population dynamics are governed by a hybrid system coupling continuous-time within-season dynamics—mortality, harvesting, and dispersal—with a discrete-time Beverton-Holt reproduction map. We derive the necessary and sufficient condition $Fr > 1$ for long-term population persistence, where $F$ encapsulates within-season survival including harvesting effects and $r$ is the intrinsic growth rate. Through bifurcation analysis, we demonstrate that marine protected areas (MPAs) significantly expand the sustainable parameter space. Using Pontryagin’s Maximum Principle, we characterize the optimal harvesting strategy as a composite Bang–Singular–Bang control. We derive an explicit state-feedback formula for the singular arc and verify its optimality via the Generalized Legendre-Clebsch condition. Numerical simulations reveal that this dynamic strategy significantly outperforms constant maximum-effort policies, yielding higher cumulative revenue while maintaining the population above the critical collapse threshold through a stable "sawtooth" trajectory. Our results highlight that modest preservation (20-30\% of habitat) allows for more intensive, profitable harvesting in open zones without risking resource extinction.

\end{abstract}

\section{Introduction}

Marine fisheries represent a vital renewable resource, supplying essential animal protein to billions of people and supporting livelihoods for hundreds of millions worldwide \cite{FAO2024}. Yet, the history of fisheries management is punctuated by severe stock collapses, exemplified by the Atlantic cod fishery off Newfoundland in the early 1990s \cite{HutchingsMyers1994,Myers1997} and the recent climate-driven declines in Pacific cod in the Gulf of Alaska \cite{McClenachan2025}. Such collapses entail profound ecological, economic, and social impacts, with recovery often proving slow and uncertain even after substantial reductions in fishing pressure \cite{Hutchings2000}. At its core, the challenge for fisheries management is bioeconomic: balancing maximization of yields and economic benefits with the preservation of population viability and ecosystem integrity.

Traditional fishery management, centered on maximum sustainable yield (MSY) principles \cite{Schaefer1954}, has been widely critiqued as inadequate for capturing the complex spatial, temporal, and multispecies dynamics of fish populations \cite{Larkin1977,Mace2001}. Contemporary approaches increasingly incorporate spatial tools, such as marine protected areas (MPAs)—particularly no-take reserves where extractive activities are prohibited \cite{Lubchenco2003,Roberts2007}. Empirical evidence and theoretical models suggest that well-designed MPAs can enhance both conservation and fishery outcomes through multiple mechanisms: providing refuges for spawning populations, protecting critical habitats, facilitating spillover of adults and larvae into adjacent fishing grounds, and reducing fishing mortality on juvenile cohorts \cite{Halpern2003, Gell2003, Roberts2001}.

The mathematical theory of renewable resource economics traces to the seminal work of \cite{Gordon1954}, who identified the "tragedy of the commons" in open-access fisheries, and \cite{Schaefer1957}, who formalized the concept of maximum sustainable yield. \cite{Clark1973} pioneered the application of optimal control theory to fishery management, demonstrating that purely economic optimization can paradoxically lead to population extinction when economic discount rates exceed biological growth rates—the so-called "Clark's rule." His monograph \cite{Clark1990} remains the definitive reference for bioeconomic modeling, establishing the theoretical foundation for much subsequent work.

Early models assumed homogeneous, non-spatial populations. Extensions incorporated age or stage structure \cite{Reed1980, Tahvonen2009}, stochastic environmental variability \cite{Reed1979, Sethi2005}, and Allee effects \cite{Courchamp2008}. Spatial heterogeneity, however, gained attention later.

\cite{SanchiricoWilen1999}  pioneered spatially explicit bioeconomic models, analyzing fishing behavior across heterogeneous patches with dispersal. They showed that spatial variation in fishing costs or stock productivity can lead to de facto reserves emerging from economic optimization, even without regulatory intervention. \cite{Neubert2003} derived the conditions under which marine reserves enhance fishery yields, identifying critical thresholds in dispersal rates and reserve sizes that determine whether reserves benefit or harm adjacent fisheries.

Empirical support for marine protected areas (MPAs) grew via meta-analyses showing consistent increases in biomass, density, size, and diversity within reserves \cite{Halpern2003, Lester2009}. Spillover to fisheries is more variable and context-dependent \cite{Harrison2012}, underscoring the value of coupled reserve-fishery models.

Recent work has extended spatial models in several directions. \cite{CostelloKaffine2010} analyzed dynamic spatial closures, showing that adaptive management of reserve boundaries can substantially improve economic outcomes. \cite{white2014value} incorporated spatial degradation and recovery dynamics, demonstrating complex trade-offs between conservation and economic objectives. \cite{BalbarMetaxas2019} reviewed ecological connectivity in MPA network design, emphasizing the importance of larval dispersal and adult movement patterns.

Many commercially important fish species exhibit pronounced seasonal patterns, with distinct periods for spawning, juvenile development, and adult feeding/migration. Classical continuous-time models implicitly assume constant vital rates and overlook these temporal structures. Discrete-time models with seasonal components better capture species like Pacific salmon \cite{Quinn2005}, Atlantic herring \cite{Overholtz2008}, and many tropical reef fish \cite{Sadovy2001}.

\cite{KarMatsuda2008} analyzed seasonal harvesting models with time delays, showing that harvesting during non-reproductive periods can sustain higher yields than year-round exploitation. \cite{Fryxell2010} demonstrated that migratory dynamics fundamentally alter optimal harvesting strategies, with implications for transboundary fisheries management. Our model builds on this tradition by explicitly separating reproduction (occurring in non-fishing season and potentially in reserves) from harvesting (occurring in designated fishing zones during limited seasons).

The Beverton-Holt recruitment function \cite{BevertonHolt1957}, which we employ, has been extensively validated for stocks exhibiting strong density-dependent recruitment at the larval or juvenile stage \cite{QuinnDeriso1999, Rose2001}. This functional form arises naturally when early life stages experience contest competition or predation that saturates at high densities, contrasting with the Ricker model appropriate for scramble competition \cite{Ricker1954}.

Pontryagin's Maximum Principle (PMP) \cite{Pontryagin1962} provides the mathematical foundation for solving dynamic optimization problems in continuous time. \cite{Clark1976} applied PMP to fishery harvesting, deriving the classic result that optimal policies often exhibit "bang-bang" control—switching between maximum and zero effort—or "singular" control, where effort follows a state-dependent feedback rule.

Our model possesses Markov structure  through its seasonal discrete-time formulation: the population state at the beginning of each harvesting season fully determines future dynamics given control choices. This enables decomposition of the infinite-horizon problem into a sequence of single-season optimizations, each solved via PMP, with the Bellman equation providing the linking condition across seasons.

In this work our contributions are:

\begin{itemize}
\item \textbf{Mathematical model}: A two-patch system coupling within-season continuous dynamics (mortality, dispersal, harvesting) with discrete-time seasonal reproduction via the Beverton-Holt function. One patch is a no-take reserve; the other allows controlled harvesting.

\item \textbf{Stability analysis}: Explicit derivation of the condition $Fr > 1$ for population persistence, where $F$ encapsulates within-season survival (including harvesting effects) and $r$ is intrinsic growth. Complete characterization of transcritical bifurcation surfaces in $(T, R, E)$ parameter space, showing how preservation zones expand stability regions.

\item \textbf{Optimal control synthesis}: Application of Pontryagin's Maximum Principle to derive optimal harvesting strategies. Proof that controls follow bang-singular-bang structure. Derivation of explicit state-feedback formula for singular arcs and verification of GLC optimality conditions.

\item \textbf{Markov decision formulation}: Demonstration that the infinite-horizon problem decomposes into seasonal subproblems linked by the Bellman equation, enabling efficient computation via value iteration or direct single-season optimization.

\item \textbf{Numerical validation}: Simulations showing that optimal dynamic controls achieve  higher revenues than maximum constant effort while maintaining sustainability. %Sensitivity analysis revealing how parameters affect both ecological stability and economic performance.
\end{itemize}

The paper is organized as follows. Section 2 formulates the two-patch seasonal model and derives the general solution for population dynamics during the harvesting season. Section 3 analyzes the non-harvested baseline system, establishing stability conditions and equilibrium properties. Section 4 extends the analysis to harvested systems, deriving the $Fr > 1$ condition and presenting comprehensive bifurcation diagrams. Section 5 formulates the bioeconomic optimization problem, proves the Markov structure, and applies Pontryagin's Maximum Principle to characterize optimal controls. Section 6 derives the explicit singular control formula and verifies optimality conditions. Section 7 presents numerical simulations comparing control strategies. Section 8 discusses ecological and economic implications, management recommendations, and model limitations. Section 9 concludes with directions for future research.

\section{Fishery Harvested Model and It's Analysis}

We assume that fishery zone is where the species arrives or passes for coming reproduction season
and the reproduction is done in the non-fishery zone and in the non-harvesting season. With these assumptions, we present by following model (system of differential equations) for adult population in the
fishery zone $x_1$ and that in the non-fishery (conserved or preserved) zone $x_2$:

For $t \in [k - T, k]$, $k = 1, 2, \dots$,

\begin{align}
\dot{x}_1 &= -c x_1 - q E x_1 + m ((1 - R) x_2 - R x_1), \\
\dot{x}_2 &= -c x_2 + m (R x_1 - (1 - R) x_2),
\label{mpaModel}
\end{align}

with $(x_1(k - T), x_2(k - T)) = ((1 - R) J_{k-1}, R J_{k-1})$.

where $T (0 \leq T \leq 1)$ is the length of harvesting season in year. Term $q E x_1$ gives the effect of harvesting in
the fishery zone with the coefficient of fishery effort $E$ and the harvesting coefficient $q$. $c$ is the natural death
rate for adult individual, $m$ is the dispersion coefficient between two habitat zones. Here $R$ $(0 \leq R < 1)$ is
the ratio of non-fishery area to the area of whole habitat zone composed with fishery and non-fishery ones.
Individuals are assumed to randomly disperse between these two zones in the habitat.

For the continuity to the other season of reproduction and maturation, we introduce $J_{k-1}$ determined by
$x(k - 1) := x_1(k - 1) + x_2(k - 1)$ as follows:

\[ J_{k-1} = \frac{r x(k - 1)}{1 + \beta x(k - 1)} \]

This gives the reproduction function for considered population by adult individuals escaped from harvesting in the $k - 1$ th year. The function is well-known Beverton-Holt type, which corresponds to the
logistic equation in continuous time model. Using this modeling with Beverton-Holt type of reproduction
function, the detail of population dynamics in reproduction season is aggregated into it, whereas it involves
the growth rate $r$ and the density effect coefficient $\beta$.

Rewriting the system (\ref{mpaModel}) in the vector form,

\[
\begin{pmatrix}
\dot{x}_1 \\
\dot{x}_2
\end{pmatrix}
= 
\begin{pmatrix}
-c - q E - m R & m (1 - R) \\
m R & -c - m (1 - R)
\end{pmatrix}
\begin{pmatrix}
x_1 \\
x_2
\end{pmatrix}
\]

i.e., $\dot{X} = A X$,

where $X = \begin{pmatrix} x_1 \\ x_2 \end{pmatrix}$, and $A = \begin{pmatrix} -c - q E - m R & m (1 - R) \\ m R & -c - m (1 - R) \end{pmatrix}$.

The solution of the system (\ref{mpaModel}) is given by,

\begin{equation} 
\label{mpaMatrixVector}
X(t) = e^{A t} X_0, 
\end{equation}

where $X_0$ can be determined by using the given initial condition. That is,

\[ X_0 = e^{-A (k - T)} X(k - T). \]

Putting this value into equation (\ref{mpaMatrixVector}), becomes

\[ X(t) = e^{-A (k - T - t)} X(k - T) \]

Thus we have

\[ X(k) = e^{A T} X(k - T), \]

i.e.,

\begin{equation} 
\label{mpaDiagonalized}
\begin{pmatrix} x_1(k) \\ x_2(k) \end{pmatrix} = e^{A T} \begin{pmatrix} 1 - R \\ R \end{pmatrix} J_{k-1} = P \begin{pmatrix} e^{\lambda_1} & 0 \\ 0 & e^{\lambda_2} \end{pmatrix} P^{-1} \begin{pmatrix} 1 - R \\ R \end{pmatrix} J_{k-1},
\end{equation}

where

\[ \lambda_i = - \frac{T}{2} \left( 2 c + m + q E + (-1)^{i+1} \sqrt{(m - q E)^2 + 4 m q E R} \right), \ (i = 1, 2) \]

is the eigenvalue of the matrix $A T$ and

\[ P = \begin{pmatrix} - m (1 - R) T & - m (1 - R) T \\ - (c + q E + m R) T - \lambda_1 & - (c + q E + m R) T - \lambda_2 \end{pmatrix} \]

is the corresponding eigenvector matrix.

From equation (\ref{mpaDiagonalized}), we have

\[ x(k) = x_1(k) + x_2(k) = F(c, q, E, m, R, T) J_{k-1} = F \frac{r x(k - 1)}{1 + \beta x(k - 1)}. \]

This is a single dimensional discrete system, where $F = e^{\lambda_1} + \frac{1}{4 m \delta} (e^{\lambda_2} - e^{\lambda_1}) (m + q E + \delta) (m - q E + \delta)$, and $\delta = \sqrt{(m - q E)^2 + 4 m q E R}$.

Let $F r = \alpha$, then we have

\begin{equation}
\label{mpaSolution}
x(k) = \frac{\alpha x(k - 1)}{1 + \beta x(k - 1)}  
\end{equation}

The equilibrium points of this discrete system are $E_1^h = 0$ and $E_2^h = \frac{\alpha - 1}{\beta}$, $\alpha > 1$. By iterations, we get the equation (\ref{mpaSolution}) as,

\[ x(k) = \frac{\alpha^k x_0}{1 + \beta [1 + \alpha + \alpha^2 + \alpha^3 + \dots + \alpha^{k-1}] x_0} = \frac{\alpha^k x_0}{1 + \beta \frac{1 - \alpha^k}{1 - \alpha} x_0} \]

As $k \to \infty$, it can be seen clearly that $x(k)$ goes to one of equilibrium states. If $0 < \alpha < 1$ then
population goes to extinction otherwise it will converge to equilibrium point $E_2^h$.

Hence it can be concluded that long term sustainable harvesting possible only when $\alpha > 1$, i.e., $F r > 1$.
Rewriting this condition, we have

\begin{equation}
\label{mpaCondition}
e^{\lambda_1} + \frac{1}{4 m \delta} (e^{\lambda_2} - e^{\lambda_1}) ((m + \delta)^2 - q^2 E^2) > 1/r 
\end{equation}

Here parameters $R$ and $E$ are controllable, so we can say that the sustainable harvesting can be done for the
suitable values of $R$ and $E$ such that the condition (\ref{mpaCondition}) satisfied. Table 1 shows the required conditions to
be hold at the extremities of $R$ and $E$.

\begin{table}[h!]
\centering
\small
\begin{tabularx}{\textwidth}{@{}l l l l X@{}}
\toprule
\textbf{Condition} & \boldmath$\delta \to$ & \boldmath$\lambda_1 \to$ & \boldmath$\lambda_2 \to$ & \boldmath$Fr > 1$ \textbf{Requirements} \\ \midrule
$R \to 0$ & $m - qE$ & $-(c+m)T$ & $-(c+qE)T$ & $E \in [0, \frac{1}{q}(\frac{\ln r}{T} - c))$ \\
$(R \neq 0)$ & & & & \textbf{OR} $T \in [0, \frac{\ln r}{c+qE})$ \\ \midrule
$R \to 1$ & $m + qE$ & $-(c+m+qE)T$ & $-cT$ & \textbf{Always holds} \\
$(R \neq 1)$ & & & & \\ \midrule
$E \to 0$ & $m$ & $-(c+m)T$ & $-cT$ & \textbf{Always holds} \\
$(E \neq 0)$ & & & & \\ \midrule
$E \to \infty$ & $\infty$ & $-\infty$ & $-[c+m(1-R)]T$ & $rRe^{-T[c+m(1-R)]} > 1$ \\
(Large $E$) & & & & \textbf{OR} $0 \leq T < \frac{\ln rR}{c+m(1-R)}$ \\
& & & & \textbf{OR} $R \gtrapprox \frac{-r + \sqrt{r^2 + 4rTme^{T(c+m)}}}{2rTm}$ \\ \bottomrule
\end{tabularx}
\caption{Condition $Fr > 1$ under extreme cases of parameters $R$ and $E$.}
\label{tab:extreme_cases}
\end{table}

\section{Non-Harvested Model and It's Analysis}

For fixed $T \in [0, 1]$, $t \in [k - T, k]$, $k = 1, 2, \dots$, we consider the following model,

\begin{align}
\label{nonHarvested}
\dot{x}_1 &= - c x_1 + m ((1 - R) x_2 - R x_1), \\
\dot{x}_2 &= - c x_2 + m (R x_1 - (1 - R) x_2),
\end{align}

with $(x_1(k - T), x_2(k - T)) = ((1 - R) J_{k-1}, R J_{k-1})$.

Note that here $T$ represents non-growth period for the population. Combining both the equations of the
system (\ref{nonHarvested}) and the initial condition, we have

\begin{equation}
\label{nonHarvested_Combined}
\dot{x}(t) = \dot{x}_1 + \dot{x}_2 = - c x_1 - c x_2 = - c x   
\end{equation}

with $x(k - T) = J_{k-1}$

Solving above initial condition problem, we get,

\[ x(k) = \exp [- c T] J_{k-1} \]

That is,

\begin{equation} 
\label{nonHarvested_Solution}
x(k) = \exp[- c T] \frac{r x(k - 1)}{1 + \beta x(k - 1)}
\end{equation}

Here $\exp[- c T] r$ is the effective intrinsic growth rate of the population considered in the model. The factor
$\exp[- c T]$ means the net death rate for the immature individual through its maturation/growth period before
the reproduction season. It is noted that $F$ (defined in previous section) become $\exp[- c T]$ if we put $E = 0$,
thus $F$ means the practical intrinsic growth rate including the harvesting effect.

The equilibrium points of this discrete system are $E_1^{nh} = 0$ and $E_2^{nh} = \frac{\exp[- c T] r - 1}{\beta}$, $r > \exp[c T]$. By
mathematical induction, we get the equation (\ref{nonHarvested_Solution}) as,

\[ x(k) = \frac{(\exp[- c T] r)^k x_0}{1 + \beta [1 + (\exp[- c T] r) + (\exp[- c T] r)^2 + (\exp[- c T] r)^3 + \dots + (\exp[- c T] r)^{k-1}] x_0}\]\[ = \frac{(\exp[- c T] r)^k x_0}{1 + \beta \frac{1 - (\exp[- c T] r)^k}{1 - (\exp[- c T] r)} x_0} \]

As $k \to \infty$, it can be seen clearly that $x_k$ goes to one of equilibrium states. If $0 < r < \exp[c T]$ then
population goes to extinction otherwise it will converge to equilibrium point $E_2^{nh}$. Hence it can be concluded
that population persist only if $r > \exp[c T](> 1)$. Note that here in non-harvested system the condition for
stability of system is dependent only on parameter $T$ (non-growth period), while in stability condition in
harvesting model is dependent on all three parameter harvesting season length ($T$), non-fishery or fishery
area ($R$ or $1 - R$) and harvesting effort $E$.

\subsection{Numerical Illustration and Bifurcation Analysis}

\begin{figure}[htbp]
    \centering
    \includegraphics[width=\textwidth]{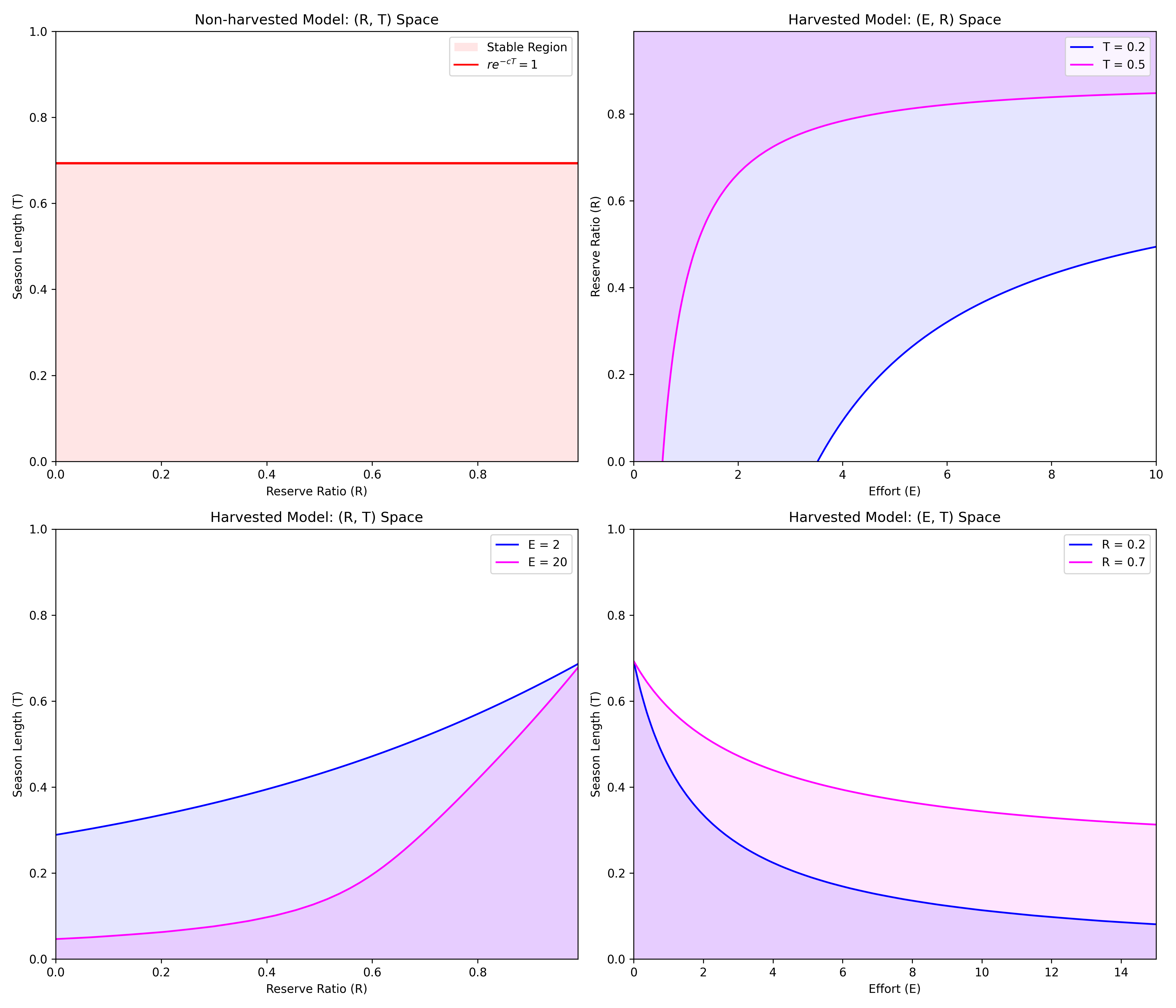}
    \caption{Bifurcation diagrams illustrating the stability boundaries for the harvested and non-harvested models. The shaded regions represent the stable parametric space where $Fr > 1$, ensuring long-term population persistence. (Top-Left): Stability in $(R, T)$ space for the non-harvested case ($E=0$). (Top-Right): Transcritical bifurcation curves in $(E, R)$ space for varying season lengths $T$. (Bottom-Left): Stability regions in $(R, T)$ space for fixed efforts $E=2$ and $E=20$. (Bottom-Right): Impact of reserve ratio $R$ on the $(E, T)$ stability boundary.}
    \label{fig:bifurcation_analysis}
\end{figure}

To analyze the parametric dependence of the system's stability, numerical simulations were conducted using a baseline set of ecologically representative parameters: $c=1,~r=2, ~m=1$, and $q=0.7$. The resulting bifurcation diagrams (Figure~\ref{fig:bifurcation_analysis}) illustrate the transcritical thresholds where the extinction equilibrium $E_1^h$ and the sustainable equilibrium $E_2^h$ exchange stability.

\begin{itemize}
    \item \textbf{Non-Harvested Stability Baseline:} In the non-harvested system ($E=0$), stability is determined solely by the condition $r \exp[-cT] > 1$. As depicted in the $(R,T)$-space, the bifurcation boundary $r \exp[-cT]=1$ appears as a horizontal line, confirming that without fishing pressure, the stability of the population depends exclusively on the non-growth period $T$ and is independent of the spatial reserve ratio $R$.
    
    \item \textbf{Reserve Fraction and Effort Capacity:} In the harvested model, the stability region is a strict subset of the non-harvested region. For a fixed harvesting season length $T$, restricting the fishery to a smaller zone (larger $R$) allows for significantly higher levels of harvesting effort $E$. Conversely, for larger fishery zones ($R \to 0$), the harvesting effort $E$ must be strictly capped at the limit $(\frac{\ln r}{T} - c)/q$ to maintain sustainability.
    
    \item \textbf{Asymptotic Stability:} A critical observation is the presence of a horizontal asymptote for the bifurcation curve $Fr=1$. If the reserve area $R$ exceeds the threshold defined by $R = (-r+\sqrt{r^{2}+4rTme^{T(c+m)}})/2rTm$, the system exhibits no upper cap for harvesting effort, implying that the protected region is large enough to sustain the global population regardless of local harvest intensity.
    
    \item \textbf{Seasonal Trade-offs:} The $(R,T)$ and $(E,T)$ parameter spaces reveal that increasing the harvesting rate $E$ in the fishery zone necessitates a compensatory reduction in the season length $T$. These simulations confirm that we cannot apply unlimited harvest rates to the system in larger fishery zones even for very brief harvesting seasons.
\end{itemize}

Collectively, these numerical results demonstrate that sustainable harvesting is strictly subjected to the interplay of parameters $T$, $R$, and $E$ satisfying the stability condition $Fr>1$. The shaded regions in Figure~\ref{fig:bifurcation_analysis} highlight the "safe" operating space for management, where biological persistence is ensured despite active extraction.

\begin{table}[h]
\centering
\small
\begin{tabularx}{\textwidth}{@{}l l X@{}}
\toprule
\textbf{Focus Area} & \textbf{Quadrant} & \textbf{Key Result and Stability Condition} \\ \midrule
Global Baseline & Top-left & Stability depends solely on $T$; $r > \exp[cT]$ must hold regardless of $R$. \\
Spatial Allocation & Top-right & Increasing $R$ (larger reserve) removes the upper cap on fishing effort $E$.  \\
Temporal Regulation & Bottom-left & Higher effort $E$ requires a shorter harvesting season $T$ to maintain $Fr > 1$. \\
System Sensitivity &  & Stability for the harvested model is a subset of the non-harvested stability region. \\ \bottomrule
\end{tabularx}
\caption{Summary of parametric dependencies and stability regions observed in Figure 1.}
\label{tab:bifurcation_summary}
\end{table}

\section{Bio-Economic Problem: The Case with Preservation Area}

In this section we formulate the bioeconomic problem; which aims, economically, to get better revenue
through fish harvesting and keep the fish population ecologically intact. To fulfil the objective, we will
determine the optimal strategy by using necessary conditions for optimality after showing existence of it.

It is common sense that more harvesting efforts generates more revenue but we may loose the population
completely in the process or it may cost us when population level is very low. Whereas keeping population
size high by doing harvesting with less effort, gives us less revenue. So the question we are interested in
is, how do we balance the harvesting effort so that we can get maximum possible revenue by harvesting
without loosing the sustainability of the population? Is there any optimal way?

This problem is named as bio-economic problem: biology demands population should be remain biologically intact, while economy demands more money. This dilemma we deal through the Pontryagin’s
Maximum Principle (PMP). The ‘main thing’ to deal these bio-economic problems is the formulation of the
problem in control theory setup, once it is done every other things are routine calculations.

[For Economical Demand] We set the objective function as:

Maximization of the revenue over admissible range of effort parameter $E$

\[ = \max_{E \in [E_{\min}, E_{\max}]} \sum_{k=1}^\infty \int_{k- T}^k (Revenue) \, dt \]

That is, maximum sum of revenue generated in every harvesting season. Since each harvesting seasons are
independent from other ones. So we can interchange the ‘max’ operator and the ‘$\sum$’ operator. In simple
words, it means that long term maximum possible revenue generated by harvesting is equal to sum of the
maximum possible revenue generated in all harvesting seasons. Then the objective function becomes

\[ = \sum_{k=1}^\infty \max_{E \in [E_{\min}(k), E_{\max}(k)]} \int_{k- T}^k (Revenue) \, dt \]

It should be noted that if we do the theoretical work to determine optimal strategy thorough PMP for
any single harvesting season. Same theoretical work be valid for other harvesting seasons, of course, with
different admissible control set as a same static admissible control set for every harvesting season may cause
the population extinction. Thus writing $[E_{\min}(k), E_{\max}(k)]$ as admissible control set for $k$-th harvesting
season will be more meaningful in problem formulation and thus we define $E_{\min} := \min_k E_{\min}(k)$ and
$E_{\max} := \max_k E_{\max}(k)$.

Also note that, it is not necessary that every element in the admissible control set $[E_{\min}, E_{\max}]$ satisfies
the asymptotic stability condition $\alpha > 1$.

In the integrand, the natural resource variable $x_1$ satisfies the proposed dynamics with given arbitrary
positive initial value. This is the given constraint, we have.

[For Ecological Demand] We need population must exist at the end of the every harvesting season.This
condition is sufficient for long-term persistence of the population.

If extinction of population happens, it will happen only in one of the harvesting seasons, since between
any two consecutive harvesting seasons population will be increased through the reproduction function $J$.
So, if we ensure the persistence by keeping certain amount population size in each harvesting season, it does
provide long-term persistence. Ensuring positive population size in every harvesting season is the guarantee
for sustainability.

\subsection{Formulation of the Optimal Control Problem}

The problem of optimal harvesting policy with preservation area $R \neq 0$ is formulated as follows:

Objective:

\begin{equation}
\max_{u \in U} \sum_{k=1}^\infty A(k) I_k, 
\end{equation} 

where $I_k = \int_{k - T}^k (P q x_1 - C) E \, dt$; net profit of harvesting in harvesting season of $k$-th year. Here, $P$ is the
per unit fish stock selling price and $C$ is the cost of harvesting. Assuming that selling price and cost are time
invariant in a particular harvesting season. Hence there is no discount rate in the integrand, while $A(k)$ is
representing the annual discount rate, that is, $(1 - \delta)^k$, where the discount rate $\delta$ assumed to be constant.
Note that with this choice of discount rate, the series $\sum^\infty_{k=1} A(k) I_k$ is convergent. Also, it can be noted that
maximization of summation is equivalent to maximization integration in summand. Harvesting period $T$ is
assumed same for every year. To formulate of the problem in optimal control theory setting, we transform
the harvesting interval $[k - T, k]$ into $[0, T]$. Thus our objective is:

\[ \max_{u \in U} \int_0^T (P q x_1 - C) E \, dt \]

Subject to the system dynamics: for $t \in [0, T] \subseteq [k - 1, k]$,

\[ \dot{x}_1 = - c x_1 - q E x_1 + m ((1 - R) x_2 - R x_1), \]

\[ \dot{x}_2 = - c x_2 + m (R x_1 - (1 - R) x_2), \]

Initial condition:

\[ (x_1(0), x_2(0)) = (x_{10}, x_{20}) = ((1 - R) J_{k-1}, R J_{k-1}), \]

= population size at the commence of

harvesting season in $k$-th year.

Terminal condition:

\[ x_1(T) = a_k > 0 \]

\[ x_2(T) = b_k > 0 \ (a_k, b_k \text{ are free, $T$ is fixed}). \]

Constraint:

\[ u(t) \in U = \{ E | E \in [E_{\min}, E_{\max}] \} \ \forall t \in [0, T] \]

Since rest periods (correspond to no harvesting) already incorporated in the model dynamics, thus we expecting at least some minimum harvesting throughout the harvesting season, that is, $E_{\min} > 0$.

\subsection{Existence and Necessary Conditions for Optimality: Bang-Bang and Singular Controls}

Prior to determine necessary conditions for optimal control, first it is always good to check the existence
of such control. For our proposed problem there exists an optimal control as sufficient conditions (refer to
[1]) are satisfied. That is, integrand of the objective function and vector field associated with state system
dynamics are continuously differentiable, $t_f = T$ is fixed with free final state, and admissible control set is
compact, convex and time invariant.

The first order necessary conditions for optimality of our formulated bioeconomic problem are determined by Pontryagin maximum principle [2]. That is, if Hamiltonian $H \equiv p^0 F_1(x, u) + p F_2(x, u)$, then

(a) Multipliers $p^0$ and $p(t)$ must not vanish simultaneously at any point of the time.

(b) Absolutely continuous co-state variable $p(t)$ follow the dynamics $\dot{p}(t) = - \frac{\partial H}{\partial x}$, and

(c) The piecewise continuous optimal control $u^*$ maximizes the Hamiltonian along the controlled state
trajectory and multipliers ($p^0, p(t)$).

It is observed that the terminal conditions are arbitrary positive real numbers, hence transversality
conditions in costate variables are $p_1(T) = 0$ and $p_2(T) = 0$. It means constant multiplier $p^0$ associated
with objective function is a positive real number, thus all the extremals for this problem are normal. We normalized the problem to unity, i.e., $p^0 = 1$, and then the Hamiltonian $H = H(p_1(t), p_2(t), x_1(t), x_2(t), E(t))$
for this problem defined as follows,

\[ H = (P q x_1 - C) E + p_1 [- c x_1 - q E x_1 + m ((1 - R) x_2 - R x_1)] + p_2 [- c x_2 + m (R x_1 - (1 - R) x_2)], \]

or equivalently,

\[ H = \langle p, f(x) \rangle + E (- C + (P q, 0) x + \langle p, g(x) \rangle), \]

where state vector $x = (x_1, x_2)^\top$, costate-vector $p = (p_1, p_2)$, and vector fields $f$ and $g$ are, respectively,
given by $(- c x_1 + m ((1 - R) x_2 - R x_1), - c x_2 + m (R x_1 - (1 - R) x_2))^\top$ and $(- q x_1, 0)^\top$. The Hamiltonian here, is a linear function of control variables thus the optimal control strategy will combination of bang-bang and singular control. Singular control corresponds to harvesting with lower than maximum effort
and higher than minimum effort, while bang-bang control represent the maximum and minimum possible
harvesting. If $E^*$ is an optimal control in the harvesting season [0, T] with corresponding state solution trajectory $x^* = (x_1^*, x_2^*)^\top$, then there exist $p_1$ and $p_2$, absolutely continuous costate variables whose dynamics
is given by Pontryagin’s maximum principle as follows,

\[ \dot{p}_1 = - \frac{\partial H}{\partial x_1} = - [P q E^* + p_1 (- c - q E^* - m R) + p_2 m R] , \ p_1(T) = 0 \]

\[ \dot{p}_2 = - \frac{\partial H}{\partial x_2} = - [p_1 m (1 - R) + p_2 (- c - m (1 - R))] , \ p_2(T) = 0 \]

i.e., 
\begin{equation}
\dot{p}(t) = - E^* (P q, 0) - \langle p(t), D f(x^*) + E^* D g(x^*) \rangle, p(T) = 0
\end{equation},

$D f$ and $D g$ are the jacobian matrices evaluated at optimal state trajectory, of the vector fields $f$ and $g$
respectively.

Also for each $t \in [0, T]$ the optimal control $E^*$ maximizes the Hamiltonian $H$ over admissible control
set $[E_{\min}, E_{\max}]$ along the optimal solution trajectory $(x_1^*, x_2^*)$ and multipliers $(p_1, p_2)$, that is,

\begin{equation} 
\label{optimal_HamiltonianCondition}
H(p_1(t), p_2(t), x_1^*(t), x_2^*(t), E^*(t)) \geq H(p_1(t), p_2(t), x_1(t), x_2(t), E(t)), 
\end{equation}

and the maximum value of $H$ is constant throughout the interval [0, T], as it is autonomous function. The
above maximum condition after simplification reduces to

\[ (- C + (P q, 0) x^* + \langle p, g(x^*) \rangle) E^* = \max_{E_{\min} \leq E \leq E_{\max}} (- C + (P q, 0) x^* + \langle p, g(x^*) \rangle) E. \]

This is equivalent to maximizing the linear function $\Phi(t) E$ over control interval $[E_{\min}, E_{\max}]$ for some
time-varying function $\Phi(t) : [0, T] \to \mathbb{R}$ known as switching function, defined as

\[ \Phi(t) = - C + (P q, 0) x^* + \langle p, g(x^*) \rangle. \]

As $H = \langle p, f(x) \rangle + E \Phi(t)$ is linear in $E$, thus the maximum of it realized for the control

\[ E^*(t) = \begin{cases} E_{\max} & \text{when } \Phi(t) > 0, \\ E_{\min} & \text{when } \Phi(t) < 0, \end{cases} \]

If switching function $\Phi(t)$ vanishes at some isolated time point $\tau$ of the zero set $L = \{ t \in [0, T] :
\Phi(t) = 0 \}$ then the control switches between $E_{\max}$ and $E_{\min}$ depending on the sign of time derivative, i.e.,
from $E_{\max}$ to $E_{\min}$ if $\dot{\Phi}(t) < 0$ and from $E_{\min}$ to $E_{\max}$ if $\dot{\Phi}(t) > 0$. However, if $\Phi(t) = 0$ holds for interval
of time, say, $[t_1, t_2]$ then the Pontryagin’s maximum principle does not provide any information about the
control $E^*(t)$. In fact in this case, all $E \in [E_{\min}, E_{\max}]$ satisfies the maximum condition (\ref{optimal_HamiltonianCondition}). The
standard procedure to get appropriate control in this situation involves solving algebraic equations resulting
from successive time derivatives of $\Phi(t)$ (as they also vanishes along $\Phi(t)$ in interval $[t_1, t_2]$) till we get
control explicitly. This type of obtained controls are called singular controls. The Generalized Legendre-
Clebsch; higher-order necessary condition for optimality need to be satisfied for these controls [3]. For
efficiently computation of time derivatives of switching function $\Phi(t)$, we will use following proposition,

Proposition 1. Let $x(t)$ is the solution of the system $\dot{x} = f(x) + E g(x)$ for control $E$ and $p(t)$ is the
solution of adjoint equation $\dot{p}(t) = - E^* (P q, 0) - \langle p(t), D (f(x^*) + E^* g(x^*)) \rangle$. Let $h(x)$ and $k(x)$ are two
smooth vector fields in $\mathbb{R}^2$ and define $\Psi(t) = (P q, 0) h(x) + \langle p, k(x) \rangle$, then the derivative of $\Psi(t)$ is given
by

\[ \dot{\Psi}(t) = (P q, 0) ((D h)(f + E g) - E k) + \langle p(t), [f + E g, k](x(t)) \rangle, \]

where $[u, v](x) = (D v(x)) u(x) - (D u(x)) v(x)$ is the Lie-bracket of the vector fields $u$ and $v$.

Proof. The proof is straight forward, applying chain rule we have,

\[ \dot{\Psi}(t) = (P q, 0) D h(x) \dot{x}(t) + \dot{p}(t) k(x) + p(t) D k(x) \dot{x}(t) \]

\[ = (P q, 0) D h \{ f + E g \} + \{ - E (P q, 0) - \langle p(t), D (f + E g) \rangle \} k + p(t) D k \{ f + E g \} \]

\[ = (P q, 0) ((D h)(f + E g) - E k) - \langle p(t), D (f + E g) \rangle k(x) + p(t) (D k)(f + E g) \]

\[ = (P q, 0) ((D h)(f + E g) - E k) + \langle p(t), [f + E g, k](x(t)) \rangle. \]

Thus using the proposition, we get $\dot{\Phi}(t) = (P q, 0) f(x) + \langle p, [f, g](x(t)) \rangle$. Since control does not appear
in the first time differentiation of the switching function, thus again applying the proposition, we have,

\[ \ddot{\Phi}(t) = (P q, 0) ((D f)(f + E g) - E [f, g](x(t))) + \langle p(t), [f + E g, [f, g]](x(t)) \rangle \]

\[ = E ((P q, 0) (D f) g - (P q, 0) [f, g](x(t) + \langle p(t), [g, [f, g]](x(t)) \rangle )) \]

\[ + (P q, 0) (D f) f + \langle p(t), [f, [f, g]](x(t)) \rangle. \]

If the coefficient of control $E$, that is, $(P q, 0) (D f) g - (P q, 0) ([f, g](x(t))) + \langle p(t), [g, [f, g]](x(t)) \rangle$ does
not vanishes in the interval $[t_1, t_2]$, then the singular control is given as follows,

\[ E^{sing} = \frac{(P q, 0) (D f) f + \langle p, [f, [f, g]](x(t)) \rangle}{(P q, 0) ([f, g](x(t))) - (P q, 0) (D f) g - \langle p(t), [g, [f, g]](x(t)) \rangle}. \]

Now, to show the existence of optimal singular control for our problem and if it exists, then writing it as
feedback function which depends upon state $x$, we require to do calculations. That goes as follows,

\[ [f, g](x) = (D g(x)) f(x) - (D f(x)) g(x) \]

\[ = \begin{pmatrix} - q & 0 \\ 0 & 0 \end{pmatrix} \begin{pmatrix} - (c + m R) x_1 + m (1 - R) x_2 \\ m R x_1 - (c + m (1 - R)) x_2 \end{pmatrix} - \begin{pmatrix} - c - m R & m (1 - R) \\ m R & - c - m (1 - R) \end{pmatrix} \begin{pmatrix} - q x_1 \\ 0 \end{pmatrix} \]

\[ = - q \begin{pmatrix} - (c + m R) x_1 + m (1 - R) x_2 \\ 0 \end{pmatrix} + q x_1 \begin{pmatrix} - (c + m R) \\ m R \end{pmatrix} = \begin{pmatrix} - q m (1 - R) x_2 \\ q m R x_1 \end{pmatrix}, \]

\begin{align*}
[f, [f, g]](x) &= \begin{pmatrix} 0 & - q m (1 - R) \\ q m R & 0 \end{pmatrix} 
\begin{pmatrix} - c x_1 + m ((1 - R) x_2 - R x_1) \\ - c x_2 + m (R x_1 - (1 - R) x_2) \end{pmatrix} \\
&\quad  -  \begin{pmatrix} - c - m R & m (1 - R) \\ m R & - c - m (1 - R) \end{pmatrix} 
\begin{pmatrix} - q m (1 - R) x_2 \\ q m R x_1 \end{pmatrix} \\[1em]
&= \begin{pmatrix} c q m (1 - R) x_2 - q m^2 R (1 - R) x_1 + q m^2 (1 - R)^2 x_2 \\ 
- c q m R x_1 + q m^2 R (1 - R) x_2 - q m^2 R^2 x_1 \end{pmatrix}\\[0.5em]
&\quad- \begin{pmatrix} c q m (1 - R) x_2 + q m^2 R (1 - R) x_2 + q m^2 R (1 - R) x_1 \\ 
- q m^2 R (1 - R) x_2 - c q m R x_1 - q m^2 R (1 - R) x_1 \end{pmatrix} \\[1em]
&= \begin{pmatrix} - 2 m^2 q R (1 - R) x_1 + q m^2 (1 - R)^2 x_2 - q m^2 R (1 - R) x_2 \\ 
2 m^2 q R (1 - R) x_2 - q m^2 R^2 x_1 + q m^2 R (1 - R) x_1 \end{pmatrix}
\end{align*}

\[ [g, [f, g]](x) = \begin{pmatrix} 0 & - q m (1 - R) \\ q m R & 0 \end{pmatrix} \begin{pmatrix} - q x_1 \\ 0 \end{pmatrix} - \begin{pmatrix} - q & 0 \\ 0 & 0 \end{pmatrix} \begin{pmatrix} - q m (1 - R) x_2 \\ q m R x_1 \end{pmatrix}\] \[= \begin{pmatrix} 0 \\ - q^2 m R x_1 \end{pmatrix} - \begin{pmatrix} q^2 m (1 - R) x_2 \\ 0 \end{pmatrix} = \begin{pmatrix} - q^2 m (1 - R) x_2 \\ - q^2 m R x_1 \end{pmatrix}, \]

\[ (P q, 0) D(f) f = P q (- c - m R, m (1 - R)) f = P q (- c - m R, m (1 - R)) \begin{pmatrix} - c x_1 + m (1 - R) x_2 - m R x_1 \\ - c x_2 - m (1 - R) x_2 + m R x_1 \end{pmatrix} \]

\[ = P q (( - c - m R) [- c x_1 + m (1 - R) x_2 - m R x_1] + m (1 - R) [- c x_2 - m (1 - R) x_2 + m R x_1]) \]

\[ = P q (c^2 x_1 - 2 c m (1 - R) x_2 + 2 c m R x_1 - m^2 R (1 - R) x_2 - m^2 (1 - R)^2 x_2 + m^2 R x_1), \]

and $(P q, 0) D(f) g = P q (- c - m R, m (1 - R)) g = P q (- c - m R, m (1 - R)) \begin{pmatrix} - q x_1 \\ 0 \end{pmatrix} = P q^2 (c + m R) x_1$

Thus the coefficient of control $E$ in $\ddot{\Phi}(t)$ is simplified as follows,

\[ \langle p, [g, [f, g]](x(t)) \rangle - (P q, 0) ([f, g](x(t))) + (P q, 0) (D f) g\]\[ = - p_1 q^2 m (1 - R) x_2 - p_2 q^2 m R x_1 + P q^2 m (1 - R) x_2 + P q^2 (c + m R) x_1\] \[= (P - p_1) m q^2 (1 - R) x_2 + (P - p_2) m q^2 R x_1 + P c q^2 x_1. \]

Since along singular optimal control, we have $\Phi(t) = 0$ and $\dot{\Phi}(t) = 0$, thus using these conditions
above expression reduces to $2 C q m (1 - R) x_2 / x_1$, which is positive always. It means that the Generalized
Legendre-Clebsch (GLC):

\[ (-1)^k \frac{\partial}{\partial E} \frac{d^{2k}}{dt^{2k}} \frac{\partial H}{\partial E} \leq 0, \ \forall t \in [t_1, t_2] \]

satisfied strictly for $k = 1$. Thus singular control is of first order [4] and all singular arcs locally maximize
the objective.

Furthermore, using conditions $\Phi(t) = 0$ and $\dot{\Phi}(t) = 0$ again along the singular extremals and then after
simplifications we have,

\[ (P q, 0) (D f) f + \langle p, [f, [f, g]](x(t)) \rangle = P q c (c + m) x_1 + 2 C m^2 R (1 - R) \left( 1 + \frac{x_2}{x_1} \right) \frac{x_2}{x_1} - 2 C m^2 (1 - R)^2 \left( 1 + \frac{x_2}{x_1} \right). \]

Hence, optimal singular control can be given as follows in terms of state variables,

\[ E^{sing}(t) = \frac{P q c (c + m) x_1 + 2 C m^2 R (1 - R) \left( 1 + \frac{x_2}{x_1} \right) \frac{x_2}{x_1} - 2 C m^2 (1 - R)^2 \left( 1 + \frac{x_2}{x_1} \right)}{2 C q m (1 - R) \frac{x_2}{x_1}} \]

\[ = \frac{P c (c + m)}{2 C m (1 - R)} \frac{x_1^2}{x_2} + \frac{m R}{q} \left( 1 + \frac{x_2}{x_1} \right) - \frac{m (1 - R)}{q} \left( 1 + \frac{x_1}{x_2} \right). \]

This optimal singular control is admissible only if it takes the values in the admissible set $[E_{\min}, E_{\max}]$ and corresponding trajectory in state space is called optimal singular subarc.

Thus, the optimal control follows a composite strategy:
$$E^*(t) = \begin{cases}
E_{max} & \text{if } E_{sing}(t) > E_{max}\\
E_{sing}(t) & \text{if } E_{min} \leq E_{sing}(t) \leq E_{max}\\
E_{min} & \text{if } E_{sing}(t) < E_{min}
\end{cases}$$

\section{Numerical Simulations}

Below are numerical simulations in support of optimal control strategies, which aim to provide best possible revenues along with ecologically intact population. Parameter values for these simulations are taken
as:

$c = 1$, $m = 1$, $q = 0.7$, $r=2$,  $R = 0.2$, $T = 0.5$, $P = 5$, $C = 1$, $\beta =0.1$, $E_{max}=10$,  and $E_{min} = 3$.

\begin{figure}[htbp]
\centering
\includegraphics[width=0.95\linewidth]{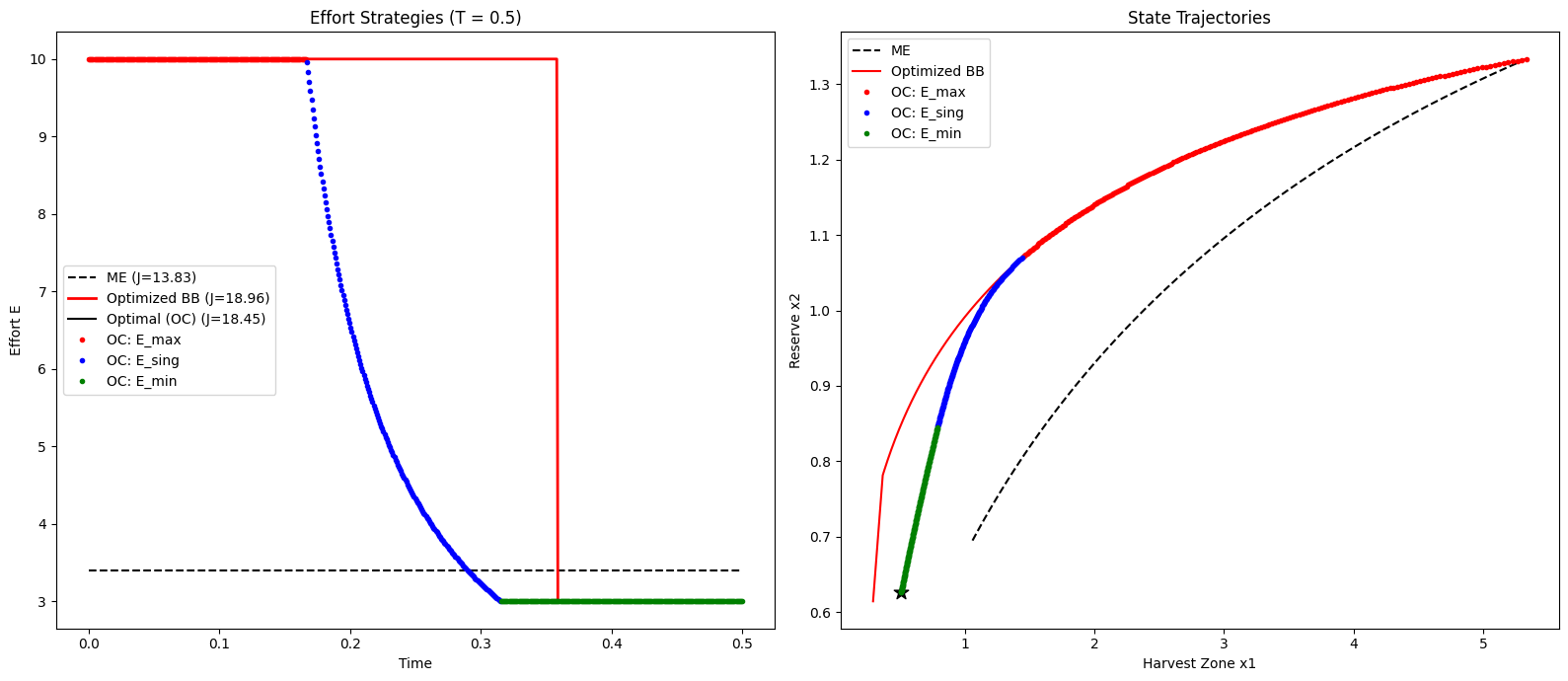}
    \caption{Comparison of harvesting effort strategies and their corresponding state trajectories. Phase (a) illustrates the effort profiles for constant control ($E=3.4$), pure bang-bang control, and the optimal combination of bang-bang and singular control ($E_{sing}$). Phase (b) presents the evolution of the population in the harvested zone ($x_1$) and reserve ($x_2$) under each strategy. The optimal control trajectory (blue) maximizes biomass extraction while ensuring the population remains within sustainable limits at the end of the season.}
\end{figure}

We evaluated the performance of three distinct harvesting strategies—Constant Effort ($ME$), Pure Bang-Bang Control, and the Optimal Control ($OC$) involving a combination of Bang-Bang and Singular arcs. Figures 2(a) and 2(b) illustrate the temporal evolution of harvesting effort and the resulting state trajectories in the $(x_1, x_2)$ phase space, respectively.

\subsection{Comparative Performance of Strategies} For a standardized single harvesting period, Figure 2(b) displays the fish population evolution under constant harvesting effort ($E = 3.4$), while the trajectories for dynamic strategies (Bang-Bang and $OC$) reflect varying intensities throughout the season. Constant Effort Strategy ($ME$): We set $E = 3.4$, which represents the largest possible constant value that satisfies the asymptotic stability condition $\alpha > 1$. Despite maintaining stability, this strategy generates significantly less revenue ($13.83$ units) compared to the optimal strategies. Optimal Control Strategy ($OC$): The $OC$ strategy utilizes the singular control $E_{\text{sing}}$ defined as an interior point within the admissible set $[E_{\text{min}}, E_{\text{max}}]$. This approach maximizes the objective function, yielding the revenue ($18.72$ and $18.45$ units) while ensuring ecological persistence. In the current parametric window ($T=0.5$), the Bang-Bang (BB) strategy shows a marginal revenue advantage (18.72 vs 18.45).

\subsection{Economic and Ecological Logic}  The reason for the suboptimal performance of the $ME$ strategy lies in its inability to adapt to fluctuating population levels. Impact of Low Population: Applying a constant $ME$ regardless of population size is costly when the stock is low, as the harvesting effort yields poor stock returns. Missed Opportunities: Conversely, when the population size is high, a constant effort fails to generate the additional revenue that would be possible by temporarily increasing intensity.

\subsection{Long-term Sustainability} As observed in the phase space trajectories (Figure 2b), the final states $(x_1(T), x_2(T))$ for all strategies are non-zero, ensuring that the population remains biologically intact. Specifically, for the $ME$ and $OC$ strategies, the terminal population size is sufficient to allow the reproduction function $J_{k-1}$ to replenish the stock for the subsequent year. This ensures that the same initial conditions can be realized annually, with the $OC$ consistently outperforming $ME$ in terms of long-term cumulative revenue.

\begin{figure}[htbp]
    \centering
    \includegraphics[width=0.95\linewidth]{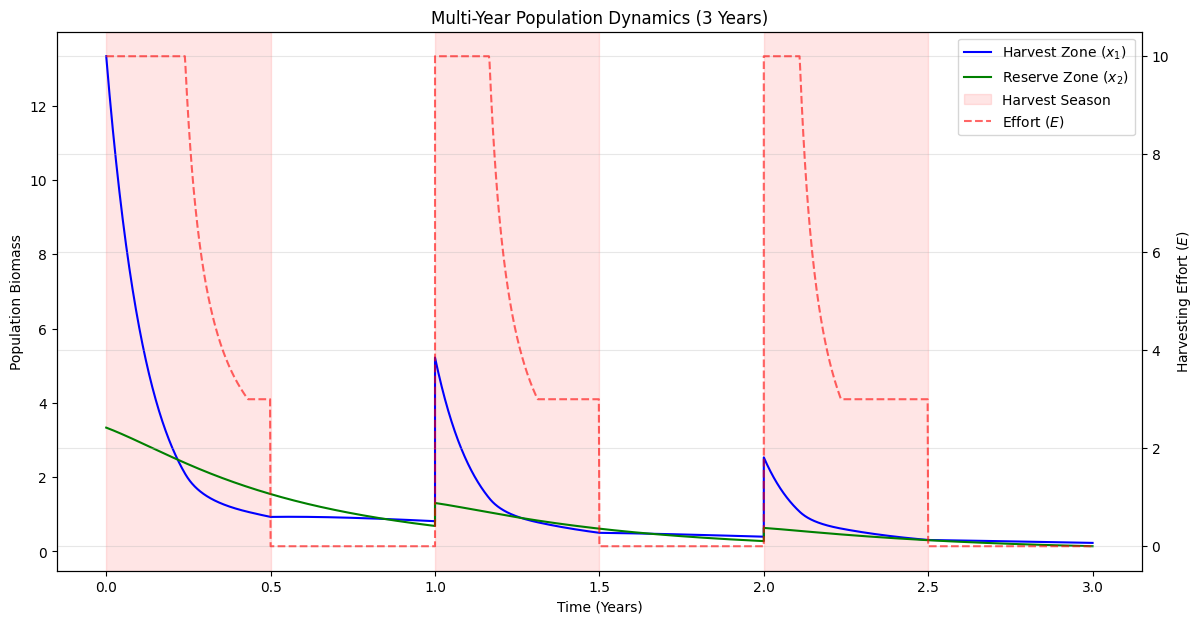}
    \caption{Multi-year bioeconomic simulation of the two-patch fishery system over a three-year horizon. The blue and green solid lines represent the population biomass in the harvested zone ($x_1$) and reserve zone ($x_2$), respectively. Shaded red regions indicate the active harvesting seasons ($T=0.5$), during which the population undergoes depletion governed by the optimal control effort $E^*$ (dashed red line). The non-shaded regions represent the recovery phase where $E=0$. Discrete vertical jumps at the start of each year ($t=1, 2, 3$) denote the inter-seasonal recruitment phase defined by the Beverton-Holt reproduction map $J$. For the chosen parameters ($r=5$), the system satisfies the sustainability condition $Fr > 1$, resulting in a stable, periodic 'sawtooth' trajectory that ensures long-term ecological persistence.}
    \label{fig:multi_year_dynamics}
\end{figure}

However, for long-term survival, the residual biomass must be sufficient to satisfy the sustainability condition $Fr > 1$. If the intrinsic growth rate $r$ is sufficiently large, this condition can be maintained despite active extraction ; if not, management must intervene by reducing the maximum harvesting pressure $E_{max}$ on the population. The multi-year dynamics illustrated in Figure \ref{fig:multi_year_dynamics} clearly demonstrate that biological persistence is not guaranteed by single-season survival alone. Instead, long-term stability requires a dynamic effort strategy that adapts to population density, ensuring that the inter-seasonal reproduction function $J_{k-1}$ can effectively replenish the stock for subsequent cycles.

\section{Discussion}
\subsection{Ecological Insights}
Our analysis reveals several key ecological findings regarding fishery management with preservation zones. Comparing bifurcation diagrams for the harvested model (Figure 1, sub-figures 2-4) with the non-harvested baseline shows that positive values of R (preservation ratio) consistently expand the parameter space permitting sustainable harvesting. This mathematical result provides theoretical support for Marine Protected Area (MPA) policies, demonstrating that even modest preservation ($R = 0.2$, or 20\% protected) significantly enhances population resilience to harvesting pressure.
Three-way parameter interaction: Unlike the non-harvested system where persistence depends solely on the balance between reproduction rate and mortality during the non-growth period ($r > e^{cT}$), the harvested system exhibits complex three-way interactions among season length ($T$), harvesting effort ($E$), and spatial structure ($R$). This complexity necessitates integrated management approaches that simultaneously consider temporal, spatial, and intensity dimensions of fishing pressure.

The dispersion rate $m$ affects stability through two competing mechanisms. Higher dispersion increases mixing between zones, which: (i) allows harvested populations to be replenished from the preserve (beneficial), but (ii) exposes protected individuals to harvesting (detrimental). Our analysis shows that the net effect depends on R, with benefits dominating when preservation zones are sufficiently large ($R > 0.15$ for our baseline parameters).
Threshold effects: The transcritical bifurcation at $Fr = 1$ represents a hard threshold separating persistence from extinction. Small decreases in reproduction rate $r$, small increases in mortality $c$, or small increases in harvesting effort $E$ can push the system across this threshold. This suggests that precautionary management is essential, as systems near the threshold exhibit little warning before collapse.

\subsection{Economic Insights}
Dynamic control substantially outperforms static strategies: Our numerical simulations (Figures 2-4) demonstrate that optimal dynamic harvesting strategies achieve revenues higher than the maximum constant-effort strategy that maintains stability. This "optimality premium" arises because dynamic strategies exploit temporal variation in population density, harvesting intensively when stocks are abundant and conservatively when stocks are depleted.

The optimal strategy exhibits a composite bang-singular-bang structure. Initially, when post-reproduction populations are high, maximum effort ($E = E_max$) is applied to capitalize on abundant stocks. As the population declines through mortality and harvesting, the control switches to a singular arc where effort tracks population density according to the feedback formula:

The optimal harvest intensity is governed by the singular control effort $E_{\text{sing}}(t)$, which integrates economic drivers and spatial biological feedback:

\begin{equation}
E_{\text{sing}}(t) = \underbrace{\frac{P c (c + m)}{2 C m (1 - R)} \frac{x_1^2}{x_2}}_{\text{Economic Incentive Term}} + \underbrace{\frac{m}{q} \left[ R \left(1 + \frac{x_2}{x_1}\right) - (1 - R) \left(1 + \frac{x_1}{x_2}\right) \right]}_{\text{Biological Stability Term}}.
\end{equation}
The \textit{Economic Incentive Term} scales effort based on the fishery's profitability ($\frac{P}{C}$) and the natural turnover rate ($c$); it identifies the harvestable surplus generated by the ``spillover'' ratio ($x_1^2/x_2$) from the reserve to the exploited zone. Conversely, the \textit{Biological Stability Term} acts as a regulatory governor; it augments effort when the reserve is robust ($x_2 > x_1$) but imposes sharp, precautionary reductions if the reserve biomass fails to sufficiently exceed the harvested zone ($x_2 \ll x_1$). Notably, the relationship $\frac{\partial E_{\text{sing}}}{\partial R} > 0$ reveals a critical bioeconomic trade-off: as the reserve fraction $R$ decreases, the optimal effort $E_{\text{sing}}$ must also be reduced. This confirms that smaller protected refugia necessitate more conservative harvesting strategies to safeguard recruitment and long-term ecological stability, effectively internalizing the precautionary principle within the feedback control.

\begin{table}[h]
\centering
\begin{tabularx}{\textwidth}{@{}l l X@{}}
\toprule
\textbf{Variable} & \textbf{Impact} & \textbf{Bioeconomic Interpretation} \\ \midrule
Price/Cost ($P/C$) & Positive ($\uparrow$) & Increased profitability justifies higher resource extraction effort. \\
Mortality ($c$) & Positive ($\uparrow$) & Higher turnover (fast-growing species) allows for more resilient, intensive harvesting. \\
Dispersal ($m$) & Positive ($\uparrow$) & Greater connectivity facilitates the transfer of biomass (spillover) from reserve to harvest zone. \\
Catchability ($q$) & Negative ($\downarrow$) & Improved efficiency means less physical effort is needed to reach the target yield. \\
Reserve Size ($R$) & Positive ($\uparrow$) & Larger protected areas provide a recruitment buffer that permits higher intensity in open zones. \\ \bottomrule
\end{tabularx}
\caption{Sensitivity analysis and bioeconomic logic of the singular control effort $E_{\text{sing}}$.}
\label{tab:bioeconomic_sensitivity}
\end{table}

\textit{Value of the preservation zone:} From a purely economic perspective, preservation zones might appear to reduce revenue by making fraction $R$ of the habitat off-limits. However, our results show that modest preservation ($R = 0.2-0.3$) actually enhances long-term economic value by: (i) expanding the range of sustainable harvesting efforts, allowing higher $E_{max}$, and (ii) providing a demographic buffer that enables more aggressive (and profitable) harvesting in the fishery zone without risking population collapse.

\textit{Shadow prices and management implications:} The costate variables $p_1(t)$ and $p_2(t)$ represent shadow prices—the marginal value of additional fish in each zone. The terminal conditions $p_1(T) = p_2(T) = 0$ indicate that populations at season's end have no direct economic value (they only matter through their contribution to next year's reproduction). This justifies aggressive late-season harvesting in the optimal strategy, but requires careful monitoring to ensure minimum viable populations persist.

\section{Conclusion}
This paper develops and analyzes an integrated bioeconomic model for fishery harvesting in spatially heterogeneous habitats with seasonal reproduction.  We derive explicit necessary and sufficient conditions for population persistence under harvesting, expressed as $Fr > 1$, where $F$ depends on mortality, dispersion, catchability, effort, season length, and reserve size. Bifurcation analysis maps stability regions in parameter space.

Using Pontryagin's Maximum Principle, we prove that optimal harvesting follows a composite bang-singular-bang strategy. We derive the singular control explicitly as a state-feedback function and verify second-order optimality conditions via the Generalized Legendre-Clebsch criterion.

Numerical simulations demonstrate that: (i) optimal dynamic strategies yield substantially higher revenues than constant-effort management, (ii) preservation zones enhance both ecological stability and economic performance, and (iii) the optimal strategy responds adaptively to population density, becoming more conservative as stocks decline.

Our results support integrated spatial-temporal management combining marine reserves with adaptive seasonal harvesting. Modest protection (20-30\% of habitat) enables more intensive and profitable harvesting in fishing zones while safeguarding against population collapse.

The model provides a theoretical foundation for ecosystem-based fisheries management, demonstrating that ecological constraints and economic objectives need not conflict when spatial refuges and dynamic control are strategically employed. Future work incorporating stochasticity, age structure, and multi-species interactions will further enhance the policy relevance of these results.

\section*{Acknowledgments}
This research problem was suggested by Professor Hiromi Seno (Tohoku University, Japan). The author gratefully acknowledges valuable discussions with Professor Seno that leads to this contribution.

\end{document}